\documentclass[12pt]{article}

\usepackage{vmargin}
\setpapersize{A4}
\setmargins{2.5cm}       
{1.5cm}                        
{16.5cm}                      
{23.42cm}                    
{10pt}                           
{1cm}                           
{0pt}                             
{2cm}                           
\usepackage{ifpdf}
\usepackage[mathscr]{eucal}
\usepackage{amsfonts}
\usepackage{amssymb}
\usepackage{color}
\usepackage {graphics}
\usepackage {epsfig}
        \usepackage{latexsym}
        \usepackage{palatino}
        \usepackage{fancybox}
\usepackage{amsmath,mathptmx,amssymb,bm}
\usepackage{amsthm}

\newenvironment{mycustom}[1]{%
  \par\vspace{3pt}%
  \noindent\textbf{\textit{#1.}}\smallskip \itshape\ 
}{\par\vspace{3pt}}

\usepackage{titlesec}

\titleformat{\paragraph}
{\normalfont\normalsize\bfseries}{\theparagraph}{1em}{}
\titlespacing*{\paragraph}
{0pt}{3.25ex plus 1ex minus .2ex}{1.5ex plus .2ex}

\usepackage{bm}

\usepackage{fancyhdr}
\usepackage{lipsum} 

\begin{document}


\title{\textbf{New non-invertible mappings and general solutions of linear wave equations with variable wave speeds}}

\author{ Rafael de la Rosa${}^{a)}$, George W. Bluman${}^{b)}$ \\ 
${}^{a)}$ Departamento de Matem\'aticas, Universidad de C\'adiz, 11510 Puerto Real,\\ C\'adiz, Spain (e-mail:rafael.delarosa@uca.es)\\
${}^{b)}$ Department of Mathematics, University of British Columbia, Vancouver, \\ British Columbia V6T 1Z2, Canada (e-mail:bluman@math.ubc.ca)}

\date{}
 
\maketitle

\thispagestyle{empty}

\begin{abstract}

We show how the symmetry-based method can be used to obtain new non-invertible equivalence mappings of linear wave equations with variable wave speeds $c(x,t)$ to linear wave equations with different variable wave speeds. Moreover, we present new non-invertible mappings of linear wave equations with variable wave speeds $c(x,t)$ to a linear wave equation with a constant wave speed. Consequently, the general solutions of these linear wave equations with variable wave speeds $c(x,t)$ are obtained.\\



\end{abstract}




\section{Introduction}

In this paper, we obtain new non-invertible mappings involving wave equations. Specifically, we present two new important results:
\begin{enumerate}
\item[1)] new non-invertible equivalence mappings of linear wave equations with variable wave speeds $c(x)$ to linear wave equations with different variable wave speeds.
\item[2)] new non-invertible mappings of linear wave equations with variable wave speeds $c(x,t)$ to a linear wave equation with a constant wave speed.
\end{enumerate}

Two systematic methods exist for constructing nonlocally related systems: the conservation law (CL)-based method and the symmetry-based method. The new results are obtained through the symmetry-based method.

Linear wave equations with variable wave speeds arise in many physical and engineering problems to study wave phenomena in an inhomogeneous continuum (see, for instance, \cite{blok2021,brek1976,mei2005}). Linear wave equations with variable wave speeds $c(x)$ appear in the analysis of small transverse vibrations of a string with variable density. Linear wave equations with variable wave speeds $c(x,t)$ can be used to model acoustic or electromagnetics waves where the wave speed changes with respect to the position and/or time due to the material or medium properties.

Different methods have been employed to determine solutions of linear wave equations with variable wave speeds. These methods include continuous symmetry and transformation methods aimed at constructing particular solutions. In \cite{blumankumei1987}, a complete point symmetry classification is given for a linear wave equation with a variable wave speed $c(x)$ and corresponding invariant solutions are determined. In \cite{blumankumeireid1988}, an algorithmic procedure is presented for determining nonlocal symmetries for PDEs in conserved form. In particular, new nonlocal symmetries of a linear wave equation with a variable wave speed $c(x)$ are found. Recently, new solutions of variable-coefficient wave equations with applications in shallow-water theory were determined in \cite{kaptsov2024,melnikov2024,pelinovsky2022}. 

The mapping of variable-coefficient linear PDEs into constant coefficient linear PDEs or nonlinear PDEs to linear PDEs is a significant area of research. Consequently, if one is able to determine solutions of these simpler PDEs, subsequently one obtains solutions of the more difficult counterparts. Several papers address this topic specifically for linear wave equations. In \cite{bluman1983}, an invertible mapping algorithm is developed to deduce whether a given linear PDE with variable coefficients can be mapped to a linear PDE with constant coefficients. This algorithm is based on mapping the infinitesimals of the Lie group of point transformations which leave invariant a given linear PDE into infinitesimals corresponding to invariance under translations. In particular, necessary and sufficient conditions are given for mapping a linear wave equation with a variable wave speed $c(x,t)$ into a linear wave equation with a constant wave speed. Moreover, the corresponding invertible mappings are given explicitly. In \cite{seymour1987},  by requiring that the asymptotic expansions employed in the theory of geometrical acoustics should terminate after a finite number of terms, the forms of the wave speed $c(x,t)$ for which any solution of a linear wave equation with this variable wave speed can be determined in terms of a corresponding solution to a linear wave equation with a constant wave speed are derived. In \cite{bluman2007}, it is proved that a linear wave equation with a variable wave speed $c(x)$ is nonlocally related to a nonlinear wave equation involving an arbitrary smooth nonlinearity.

The symmetry-based method for systematically obtaining nonlocally related PDE systems and finding nonlocal symmetries for a given PDE system with two independent variables from any admitted point symmetry is introduced in \cite{blumanyang2013,yang2013}. By taking into account the canonical coordinates for an admitted point symmetry and considering the translated canonical coordinate as a new dependent variable, a nonlocally related PDE system can be constructed. 

In \cite{blumanrafamsml2021a,blumanrafamsml2021b}, we show that in the symmetry-based method each derivative of the translated canonical coordinate is a differential invariant corresponding to a group of translations. As a result, the symmetry-based method can be naturally extended to ordinary differential equations (ODEs) and to PDEs with any number of independent variables. Moreover, in \cite{blumanrafamsml2021a}, we establish a connection between the symmetry-based method for determining nonlocally related systems for PDEs and Lie's reduction of order algorithm for ODEs. This connection is further explored in \cite{blumanrafa2025}. Here, we revisit Lie's reduction of order algorithm for ODEs from a new perspective which shows that the relationship between a given ODE and its reduced ODE is nonlocal. More significantly, we show that the symmetry-based method for determining nonlocally related PDE systems of a given PDE with any number of independent variables is a direct and natural generalization to PDEs of Lie's reduction of order algorithm for ODEs.

An application of the symmetry-based method to linear PDE systems is shown in \cite{blumanzuhal}. In particular, a systematic method is presented for determining non-invertible mappings of linear PDEs to linear PDEs. Consequently, new exact solutions of a given linear PDE can be obtained from known solutions of a nonlocally equivalent linear PDE. Note that given a linear PDE system, applying the theory presented in \cite{blumanrafamsml2021a,blumanyang2013,yang2013} is essentially ineffective, since, in general, the resulting nonlocally related PDE system is nonlinear.

A noteworthy situation arises for variational linear PDE systems (self-adjoint linear PDE systems). Here, one can directly obtain three distinct nonlocally related linear PDE systems from an admitted point symmetry. This follows from noticing that a solution of the linear PDE system yields a nonlocally related PDE system from both the CL-based and symmetry-based methods, which are not necessarily invertibly equivalent.

The CL-based method was originally presented in \cite{blumankumei1987} and subsequently broadened in \cite{blumankumeireid1988}, where a systematic procedure was proposed for deriving nonlocally related PDE systems and consequent nonlocal symmetries (known as potential symmetries) from any admitted CL of a given PDE system. This method was later generalized to PDE systems with three or more independent variables and PDE systems admitting more than one CL. For additional details, see \cite{ancobluman1997,ancothe2005,bluman2007,blumanchevianco,bluman2006,blumankumei}.

In \cite{blumanrafamsml2021b}, we establish direct connections between the symmetry-based and CL-based methods. More importantly, for PDE systems with two independent variables, we show that the CL-based method is a special case of the symmetry-based method.  

This paper is organized as follows. In Section \ref{sec:equivmap}, we present equivalence mappings of linear wave equations with variable wave speeds $c(x)$ to linear wave equations with different variable wave speeds. In Section \ref{sec:invmap}, we focus on invertible mappings. Here we present invertible equivalence mappings for linear wave equations with variable wave speeds previously obtained in \cite{blumanrafamsml2021a,blumanrafamsml2021b} and a new one. In Section \ref{sec:noninvmap}, we focus on non-invertible equivalence mappings obtained by the symmetry-based method. Here we include such an equivalence mapping obtained in \cite{blumanrafamsml2021b} as well as a new one.

In Section \ref{sec:vartocons}, we exhibit mappings of linear wave equations with variable wave speeds $c(x,t)$ to a linear wave equation with a constant wave speed. Analogously, in Section \ref{sec:invmapcte}, we focus on invertible mappings and review previous results presented in \cite{bluman1983}. In Section \ref{sec:noninvmapcte}, we present new non-invertible mappings of linear wave equations with variable wave speeds $c(x,t)$ to a linear wave equation with a constant wave speed. Moreover, we determine the general solutions of these linear wave equations with variable wave speeds $c(x,t)$.

In Section \ref{sec:genresult}, we show how the non-invertible equivalence mappings presented in Section \ref{sec:noninvmap} are obtained through the symmetry-based method. 

Finally, we present concluding remarks.

\section{Equivalence mappings for linear wave equations with variable wave speeds}\label{sec:equivmap}

Consider the linear wave equation 
\begin{equation}
\label{waveeqcx}
\dfrac{u_{tt}}{c^2(x)}=u_{xx},
\end{equation}
with a variable wave speed $c(x)$. Let $u=f(x,t)$ be any solution of the linear wave equation (\ref{waveeqcx}). Since the linear operator 
\begin{equation}\label{linerop} L=\dfrac{1}{c^2(x)}\dfrac{\partial^2}{\partial t^2}-\dfrac{\partial^2}{\partial x^2},\end{equation} 
is self-adjoint, $f(x,t)\frac{\partial}{\partial u}$ is a point symmetry of the linear wave equation (\ref{waveeqcx}) and $f(x,t)$ is a multiplier for a conservation law of the linear wave equation (\ref{waveeqcx}). It follows that a solution of the linear wave equation (\ref{waveeqcx}) yields directly a nonlocally related linear PDE system for both the symmetry-based and CL-based methods.

In \cite{blumanrafamsml2021b}, we show that for the linear wave equation (\ref{waveeqcx}), with $c(x)$ arbitrary, the nonlocally related linear PDE systems obtained by both methods for the particular solutions $f(x,t)=1$ and $f(x,t)=x$, are one-to-one equivalent. Moreover, no new equivalence mappings of the linear wave equation (\ref{waveeqcx}) can be obtained for the solution $f(x,t)=1$. 

Equivalence mappings for linear wave equations with variable wave speeds that are nonlocally related to the linear wave equation (\ref{waveeqcx}), with $c(x)$ arbitrary, are systematically constructed from some particular solutions of the linear wave equation (\ref{waveeqcx}) through the symmetry-based method in \cite{blumanrafamsml2021a,blumanrafamsml2021b}. Here, in addition to reviewing previous results derived in \cite{blumanrafamsml2021a,blumanrafamsml2021b}, we present a new equivalence mapping for a linear wave equation with a variable wave speed nonlocally related to the linear wave equation (\ref{waveeqcx}) for any wave speed $c(x)$.

\subsection{Invertible mappings}\label{sec:invmap}

As mentioned in the Introduction, valuable nonlocally related PDE systems can be obtained by using the CL-based method. In \cite{blumanrafamsml2021b}, we determine the corresponding nonlocally related equivalent linear PDE systems for the linear wave equation (\ref{waveeqcx}) with a variable wave speed $c(x)$ arising from the CL-based method when one considers its particular solutions $u=1$, $u=x$ and $u=t$, respectively. However, if one applies the CL-based method to seek non-invertible equivalence mappings of the linear wave equation (\ref{waveeqcx}) with a variable wave speed $c(x)$ to a linear wave equation with a different variable wave speed, the resulting equivalence mappings turn out to be invertible. Here, in addition to reviewing previous results obtained in \cite{blumanrafamsml2021a,blumanrafamsml2021b}, we include a new linear wave equation with a variable wave speed that is locally related to the linear wave equation (\ref{waveeqcx}) with a variable wave speed $c(x)$, and obtain the corresponding invertible mapping. 

\medskip

\begin{mycustom}{Invertible equivalence mappings of the linear wave equation (\ref{waveeqcx}) with a variable wave speed $\boldsymbol{c(x)}$ to a linear wave equation with a different variable wave speed}\label{resultCLcarb}

The linear wave equation (\ref{waveeqcx}) with a variable wave speed $c(x)$ is locally related to: 
\begin{itemize}
\item[(1)] a linear wave equation with a variable wave speed $\tilde{c}(x)=x^2 c\left(-\dfrac{1}{x}\right)$. The mapping
\begin{equation}
\label{mapcxsolxCL}
x=-\frac{1}{X}, \qquad t=t, \qquad u=-\frac{1}{X} \sigma\left(X, t \right),
\end{equation}
transforms the linear wave equation (\ref{waveeqcx}) to the linear wave equation given by
\begin{equation}\label{finalrelatedwaveCLsystemcxf=x}\dfrac{\sigma_{tt}}{X^4 c^2\left(-\dfrac{1}{X}\right)}=\sigma_{XX}.\end{equation}
This result appears in \cite{blumanrafamsml2021a}.


\item[(2)] a linear wave equation with a variable wave speed $\tilde{c}(x,t)=t^{-2} c\left(x\right)$. The mapping
\begin{equation}
\label{mapcxsoltCL}
x=x, \qquad t=-\dfrac{1}{T}, \qquad  u=-\dfrac{1}{T} \sigma\left(x, T \right),
\end{equation}
transforms the linear wave equation (\ref{waveeqcx}) to the linear wave equation given by 
\begin{equation}\label{finalrelatedwaveCLsystemcxf=t}\dfrac{\sigma_{TT}}{T^{-4} c^2\left(x\right)}=\sigma_{xx}.\end{equation}
This result appears in \cite{blumanrafamsml2021b}.

\item[(3)] a linear wave equation with a variable wave speed $\tilde{c}(x,t)=x^2 t^{-2} c\left(-\dfrac{1}{x}\right)$. The mapping
\begin{equation}
\label{mapcxsolx2}
x=-\frac{1}{X}, \qquad t=-\frac{1}{T},   \qquad u=\frac{1}{X T} \sigma\left(X,T\right),
\end{equation}
transforms the linear wave equation (\ref{waveeqcx}) to the linear wave equation given by
\begin{equation}\label{finalrelatedwaveDIsystemcxf=x2}\dfrac{\sigma_{TT}}{X^4 T^{-4} c^2\left(-\dfrac{1}{X}\right)}=\sigma_{XX}.\end{equation} 
This is a new result.
\end{itemize}
\end{mycustom}

\subsection{Non-invertible mappings}\label{sec:noninvmap}

In this section, we present new non-invertible equivalence mappings of the linear wave equation (\ref{waveeqcx}) with a variable wave speed $c(x)$ to a linear wave equation with a different variable wave speed through the symmetry-based method. In \cite{blumanrafamsml2021a,blumanrafamsml2021b}, we determine the corresponding nonlocally related equivalent linear systems for the linear wave equation (\ref{waveeqcx}) with a variable wave speed $c(x)$ arising from the symmetry-based method when one considers its particular solutions $u=1$, $u=x$ and $u=t$, respectively. Furthermore, the symmetry-based method could yield a non-invertible equivalence mapping of the linear wave equation (\ref{waveeqcx}) with a variable wave speed $c(x)$ to a linear wave equation with a different variable wave speed from any particular solution of the linear wave equation (\ref{waveeqcx}). 

In \cite{blumanrafamsml2021b}, we show that if one chooses the solutions $u=1$ and $u=x$ of the linear wave equation (\ref{waveeqcx}), the nonlocally related linear PDE systems obtained by both the CL-based and symmetry-based methods are one-to-one equivalent. Consequently, we focus on the particular solution $u=t$. Here, in addition to reviewing previous results presented in \cite{blumanrafamsml2021b}, we include a new linear wave equation with a variable wave speed that is nonlocally related to the linear wave equation (\ref{waveeqcx}) with a variable wave speed $c(x)$ through the symmetry-based method, and obtain the corresponding non-invertible mappings. 

\medskip

\begin{mycustom}{Non-invertible equivalence mappings of the linear wave equation (\ref{waveeqcx}) with a variable wave speed $\boldsymbol{c(x)}$ to a linear wave equation with a different variable wave speed}\label{resultDIcarb}

The linear wave equation (\ref{waveeqcx}) with a variable wave speed $c(x)$ is nonlocally related through the symmetry-based method to:
\begin{itemize}
\item[(1)] a linear wave equation with a variable wave speed $\tilde{c}(x,t)=t^{-4/3} c\left(x\right)$. The mapping
\begin{equation}
\label{mapcxsolt2}
x=x, \qquad t=3 T^{-1/3}, \qquad  u=T^{-1/3} \int T^{-5/3} B\left(x,T\right) \, dT,
\end{equation}
transforms the linear wave equation (\ref{waveeqcx}) to the linear PDE given by
\begin{equation}\label{intermapcxsolt2} \frac{T B_T- B}{c^2\left(x\right)} = \int  T^{-5/3} B_{xx} \, dT.\end{equation}
After differentiating PDE (\ref{intermapcxsolt2}) with respect to $T$, we obtain the linear wave equation 
\begin{equation}\label{finalrelatedwaveDIsystemcxf=t2}\dfrac{B_{TT}}{T^{-8/3} c^2\left(x\right)}=B_{xx}.\end{equation} 
This result appears in \cite{blumanrafamsml2021b}.
\item[(2)] a linear wave equation with a variable wave speed $\tilde{c}(x,t)=t^{-2/3} c\left(x\right)$. The mapping
\begin{equation}
\label{mapcxsolt1}
x=x, \qquad t=3 T^{1/3}, \qquad  u=T^{1/3} \int T^{-4/3} B\left(x,T\right) \, dT,
\end{equation}
transforms the linear wave equation (\ref{waveeqcx}) to the linear PDE given by
\begin{equation}\label{intermapcxsolt1} \frac{B_T}{c^2\left(x\right)} = \int  T^{-4/3} B_{xx} \, dT.\end{equation}
After differentiating PDE (\ref{intermapcxsolt1}) with respect to $T$, we obtain the linear wave equation 
\begin{equation}\label{finalrelatedwaveDIsystemcxf=t1}\dfrac{B_{TT}}{T^{-4/3} c^2\left(x\right)}=B_{xx}.\end{equation} 
This is a new result.
\end{itemize}
\end{mycustom}
The proofs of these results are presented in Section \ref{sec:genresult}.


\section{Mappings of linear wave equations with variable wave speeds to a linear wave equation with a constant wave speed}\label{sec:vartocons}

Consider the linear wave equation
\begin{equation}
\label{waveeq}
\dfrac{u_{tt}}{c^2(x,t)}=u_{xx},
\end{equation} 
with a variable wave speed $c(x,t)$, and the linear wave equation with a constant wave speed
\begin{equation}
\label{ctecoeffpde}
V_{\xi \eta}=0.
\end{equation}

We now focus on the construction of mappings of the linear wave equation (\ref{waveeq}) with a variable wave speed $c(x,t)$ to the linear wave equation with a constant wave speed (\ref{ctecoeffpde}). In \cite{bluman1983}, all wave speeds $c(x,t)$, as well as the corresponding mappings, for which the linear wave equation (\ref{waveeq}) can be invertibly mapped to the linear wave equation with a constant wave speed (\ref{ctecoeffpde}) are determined. For a new class of wave speeds $c(x,t)$ we determine non-invertible mappings of the linear wave equation (\ref{waveeq}) to the linear wave equation with a constant wave speed (\ref{ctecoeffpde}) by combining the results proved in \cite{bluman1983} with the results obtained in Section \ref{sec:noninvmap}. Consequently, we obtain the general solutions of these linear wave equations with the new class of variable wave speeds $c(x,t)$. 

\subsection{Invertible mappings}\label{sec:invmapcte}

We now review previous results proved in \cite{bluman1983} concerning the construction of invertible mappings of the linear wave equation (\ref{waveeq}) for all allowable wave speeds $c(x,t)$ to the linear wave equation with a constant wave speed (\ref{ctecoeffpde}). We present these results modulo scalings and translations of $x$ and $t$ to remove redundant constants.  

\medskip

\begin{mycustom}{Invertible mappings of the linear wave equation (\ref{waveeq}) to the linear wave equation with a constant wave speed (\ref{ctecoeffpde}) for all allowable wave speeds $\boldsymbol{c(x,t)}$ }

The linear wave equation (\ref{waveeq}) with a variable wave speed $c(x,t)$ can be invertibly mapped to the linear wave equation with a constant wave speed (\ref{ctecoeffpde}) if and only if the wave speed is of the form
\begin{equation}
\label{wavespeed}
c(x,t)=\frac{x^2-\Delta}{t^2-\Delta},
\end{equation}
where $\Delta$ is an arbitrary constant. The cases $\Delta=0$, $\Delta>0$ and $\Delta<0$, are considered separately.
\begin{enumerate}
\item If $\Delta=0$, the mapping is given by
$$\xi=\frac{1}{x}+\frac{1}{t}, \qquad \eta=\frac{1}{x}-\frac{1}{t}, \qquad V(\xi,\eta)=\frac{u(x,t)}{x \, t}.$$
\item If $\Delta=\rho^2>0$, the mapping is given by
$$\begin{array}{rcl}\xi &=&\log \left|\dfrac{x-\rho}{x+\rho} \right|^{\rho/2}-\log \left|\dfrac{t-\rho}{t+\rho} \right|^{\rho/2}, \vspace*{0.2cm} \\
 \eta&=&\log \left|\dfrac{x-\rho}{x+\rho} \right|^{\rho/2}+\log \left|\dfrac{t-\rho}{t+\rho} \right|^{\rho/2}, \vspace*{0.2cm} \\ 
 V(\xi,\eta)&=&\big( \left(x^2-\rho^2\right)\left(t^2-\rho^2 \right) \big)^{-1/2}\,u(x,t).
 \end{array}$$
\item If $\Delta=-\rho^2<0$, the mapping is given by
$$\begin{array}{rcl}\xi &=&\rho \arctan \left( \dfrac{x}{\rho} \right) -\rho \arctan \left( \dfrac{t}{\rho} \right), \vspace*{0.2cm} \\
 \eta&=& \rho\arctan \left( \dfrac{x}{ \rho} \right) +\rho \arctan \left( \dfrac{t}{\rho} \right), \vspace*{0.2cm} \\ 
 V(\xi,\eta)&=&\big( \left(x^2+\rho^2\right)\left(t^2+\rho^2 \right) \big)^{-1/2} \,u(x,t).
 \end{array}$$
\end{enumerate}
Moreover, in the special case $c(x,t)=c(x)$, the linear wave equation (\ref{waveeqcx}) can be invertibly mapped to the linear wave equation with a constant wave speed (\ref{ctecoeffpde}) if and only if $c(x)=x^2$, i.e., the linear wave equation (\ref{waveeqcx}) takes the form
\begin{equation}
\label{waveeqquadcx}
\dfrac{u_{tt}}{x^4}= u_{xx}.
\end{equation}
Here the mapping is
\begin{equation}\label{mapcx} \xi=\frac{1}{x}+t, \qquad \eta=\frac{1}{x}-t, \qquad V(\xi,\eta)=\frac{u(x,t)}{x}.\end{equation}
Hence the general solution of the linear wave equation (\ref{waveeqquadcx}) is given by
$$u(x,t)=x \left( F\left(\frac{1}{x}+t\right)+G\left(\frac{1}{x}-t\right) \right),$$
with $F$ and $G$ arbitrary functions of their respective arguments at least of class ${\cal C}^2$.
\end{mycustom}

\subsection{Non-invertible mappings}\label{sec:noninvmapcte}

The results reviewed in Section \ref{sec:invmapcte} serve as a starting point for determining non-invertible mappings of the linear wave equation (\ref{waveeq}) to the linear wave equation with a constant wave speed (\ref{ctecoeffpde}), for a new class of variable wave speeds $c(x,t)$. In turn, one obtains  the general solutions of this new class of linear wave equations with variable wave speeds $c(x,t)$. Here, we focus on the linear wave equation (\ref{waveeqcx}) with $c(x)=x^2$, i.e., the linear wave equation (\ref{waveeqquadcx}). The following new important result is obtained.

\medskip

\begin{mycustom}{Non-invertible mappings of the linear wave equation (\ref{waveeq}) for a new class of wave speeds $\boldsymbol{c(x,t)}$ to the linear wave equation with a constant wave speed (\ref{ctecoeffpde})}

The linear wave equation with a constant wave speed (\ref{ctecoeffpde}) is invertibly related to the linear wave equation (\ref{waveeqquadcx}). Moreover, the linear wave equation (\ref{waveeqquadcx}) is nonlocally related to (1) a linear wave equation with a variable wave speed $\tilde{c}(x,t)=t^{-4/3} x^2$; (2) a linear wave equation with a variable wave speed $\tilde{c}(x,t)=t^{-2/3} x^2$. Consequently, it follows that the linear wave equation with a constant wave speed (\ref{ctecoeffpde}) is nonlocally related to:
\begin{itemize}
\item[(1)] a linear wave equation with a variable wave speed $\tilde{c}(x,t)=t^{-4/3} x^2$. The mapping 
\begin{equation}\label{mapcxsolt1DIsimplest}
\left\{ \begin{array}{rcl}
\xi &=& \dfrac{1}{x}+3 T^{-1/3}, \vspace*{0.15cm}\\
\eta &=&  \dfrac{1}{x}-3 T^{-1/3}, \vspace*{0.15cm}\\
V &=& \dfrac{T^{-1/3}}{x} \displaystyle \int T^{-5/3} B(x,T) \, dT,
\end{array} \right.
\end{equation}
arising from the combination of mappings (\ref{mapcxsolt2}) and (\ref{mapcx}), transforms the linear wave equation (\ref{ctecoeffpde}) to the linear PDE given by
\begin{equation}\label{intermapcxsolt2quad} \frac{T B_T- B}{x^4} = \int  T^{-5/3} B_{xx} \, dT.\end{equation}
After differentiating PDE (\ref{intermapcxsolt2quad}) with respect to $T$, we obtain the linear wave equation 
\begin{equation}\label{finalrelatedwaveDIsystemcxf=t2q}
\dfrac{B_{TT}}{T^{-8/3} x^4}=B_{xx}.
\end{equation}
Hence, the general solution of the linear wave equation (\ref{finalrelatedwaveDIsystemcxf=t2q}) is given by
$$\begin{array}{rcl} B(x,T) &=& x \, T \Bigg( F \left(\dfrac{1}{x}+3 T^{-1/3} \right)+G\left(\dfrac{1}{x}-3 T^{-1/3} \right) 
\vspace*{0.15cm} \\ & &  \qquad + 3 T^{-1/3} \Bigg( G'\left( \dfrac{1}{x}-3 T^{-1/3} \right)- F'\left(\dfrac{1}{x}+3 T^{-1/3}\right) \Bigg) \Bigg), \end{array}$$
with $F$ and $G$ arbitrary functions of their respective arguments at least of class ${\cal C}^3$.

\item[(2)] a linear wave equation with a variable wave speed $\tilde{c}(x,t)=t^{-2/3} x^2$. The mapping
\begin{equation}\label{mapcxsolt2DIsimplest}
\left\{ \begin{array}{rcl}
\xi &=& \dfrac{1}{x}+3 T^{1/3}, \vspace*{0.15cm}\\
\eta &=&  \dfrac{1}{x}-3 T^{1/3}, \vspace*{0.15cm}\\
V &=& \dfrac{T^{1/3}}{x} \displaystyle \int T^{-4/3} B(x,T) \, dT,
\end{array} \right.
\end{equation}
arising from the combination of mappings (\ref{mapcxsolt1}) and (\ref{mapcx}), transforms the linear wave equation (\ref{ctecoeffpde}) to the linear PDE given by
\begin{equation}\label{intermapcxsolt1quad} \frac{B_T}{x^4} = \int  T^{-4/3} B_{xx} \, dT.\end{equation}
After differentiating PDE (\ref{intermapcxsolt1quad}) with respect to $T$, we obtain the linear wave equation 
\begin{equation}\label{finalrelatedwaveDIsystemcxf=t1q}
\dfrac{B_{TT}}{T^{-4/3} x^4}=B_{xx}.
\end{equation}
Hence, the general solution of the linear wave equation (\ref{finalrelatedwaveDIsystemcxf=t1q}) is given by
$$\begin{array}{rcl} B(x,T) &=& x \Bigg( F \left(\dfrac{1}{x}+3 T^{1/3}\right)+G\left(\dfrac{1}{x}-3 T^{1/3} \right) 
\vspace*{0.15cm} \\ & &  \qquad + 3 T^{1/3} \Bigg( G'\left( \dfrac{1}{x}-3 T^{1/3} \right)- F'\left(\dfrac{1}{x}+3 T^{1/3}\right) \Bigg) \Bigg), \end{array}$$
with $F$ and $G$ arbitrary functions of their respective arguments at least of class ${\cal C}^3$.
\end{itemize}
\end{mycustom}

\section{Proofs}\label{sec:genresult}

In this section, we show how the symmetry-based method, for the solution $u=t$ of the linear wave equation (\ref{waveeqcx}), yields the non-invertible equivalence mappings of the linear wave equation (\ref{waveeqcx}) with a variable wave speed $c(x)$ to linear wave equations with different variable wave speeds, presented in Section \ref{sec:noninvmap}.

Let $u=f(x,t)$ be any solution of the linear wave equation (\ref{waveeq}) and let $U=\frac{u}{f(x,t)}$. In $(x,t,U)$-coordinates, the linear wave equation (\ref{waveeq}) becomes the invertibly related linear PDE
\begin{equation}
\label{relatedwaveeq}
\dfrac{2 f_t U_t+f U_{tt}}{c^2(x,t)}=2 f_x U_x+f U_{xx}.
\end{equation}
From the form of the linear PDE (\ref{relatedwaveeq}), we can naturally introduce as auxiliary dependent variables the differential invariants $\alpha=U_x$, $\beta=U_t$. Then the linear PDE (\ref{relatedwaveeq}) becomes the nonlocally related linear PDE system 
\begin{equation}
\label{relatedwavesystem}
\left\{ \begin{array}{rcl}
\beta_x&=&\alpha_t, \vspace*{0.1cm}\\ 
\dfrac{2 f_t \beta+f \beta_{t}}{c^2(x,t)}&=&2 f_x \alpha+f \alpha_{x}.
\end{array}\right.
\end{equation}
Now, let $\alpha=f^{-2}\varphi(x,t)$. Then the nonlocally related linear PDE system (\ref{relatedwavesystem}) becomes 
\begin{equation}
\label{relatedwavesystemcxt}
\left\{ \begin{array}{rcl}
\beta_x&=&(f^{-2}\varphi)_t, \vspace*{0.1cm}\\ 
\dfrac{\left(f^2\beta\right)_t}{c^2(x,t)}&=&\varphi_x.
\end{array}\right.
\end{equation}


Consider the special case $c(x,t)=c(x)$. Here, the nonlocally related linear PDE system (\ref{relatedwavesystemcxt}) takes the form
\begin{equation}
\label{relatedwavesystem2}
\left\{ \begin{array}{rcl}
\beta_x&=&(f^{-2}\varphi)_t, \vspace*{0.1cm}\\ 
\dfrac{\left(f^2\beta\right)_t}{c^2(x)}&=&\varphi_x.
\end{array}\right.
\end{equation}

We now focus on the solution $u=f(x,t)=t$ of the linear wave equation (\ref{waveeqcx}). After multiplying the second equation of the linear PDE system (\ref{relatedwavesystem2}) by $t^{-2}$, the linear PDE system (\ref{relatedwavesystem2}) becomes
\begin{equation}
\label{relatedwavesystem3cxf=t}
\left\{ \begin{array}{rcl}
\displaystyle \beta_x &=& \left( t^{-2} \varphi \right)_t, \vspace*{0.15cm} \\ 
\displaystyle \frac{t^{-2} \left( t^2\beta \right)_{t}}{c^2(x)} &=& \left( t^{-2} \varphi \right)_x.
\end{array}\right.
\end{equation}
Excluding $\varphi$ through cross-differentiation, the linear PDE system (\ref{relatedwavesystem3cxf=t}) yields the scalar linear PDE 
\begin{equation}
\label{scalarwaveeqcxf=t}
\dfrac{\left(t^{-2}\left(t^2\beta\right)_t\right)_t}{c^2(x)}=\beta_{xx}.
\end{equation}
Let $\beta=\left( \gamma_1 t^{-2}+\gamma_2 t \right) B(x,t)$, with $\gamma_1$ and $\gamma_2$ arbitrary constants not simultaneously zero. Then the linear PDE (\ref{scalarwaveeqcxf=t}) becomes
\begin{equation}
\label{scalarwaveeqcxf=t2} \displaystyle \frac{1}{c^2(x)}\left(\frac{2 \left( 2 t^3 \gamma_2-\gamma_1 \right) }{t \left( t^3 \gamma_2+\gamma_1 \right)}B_t+ B_{tt} \right)= B_{xx}.
\end{equation}
Now, we distinguish the cases $\gamma_2 \neq 0$, $\gamma_2 = 0$. 

If $\gamma_2 \neq 0$, without loss of generality $\gamma_2=1$, $\gamma_1=0$. Let $T=27 t^{-3}$. Then, the linear PDE (\ref{scalarwaveeqcxf=t2}) becomes the linear wave equation (\ref{finalrelatedwaveDIsystemcxf=t2}). Thus, through the symmetry-based method, we have shown that the linear wave equation (\ref{waveeqcx}) with a variable wave speed $c(x)$ is nonlocally equivalent to a linear wave equation with a variable wave speed $\tilde{c}(x,t)=t^{-4/3} c\left(x\right)$.

If $\gamma_2=0$, without loss of generality $\gamma_1=1$. Let $T=\dfrac{t^3}{27}$. Then, the linear PDE (\ref{scalarwaveeqcxf=t2}) becomes the linear wave equation (\ref{finalrelatedwaveDIsystemcxf=t1}). Thus, through the symmetry-based method, we have shown that the linear wave equation (\ref{waveeqcx}) with a variable wave speed $c(x)$ is nonlocally equivalent to a linear wave equation with a variable wave speed $\tilde{c}(x,t)=t^{-2/3} c\left(x\right)$. This is a new result.

Consequently, we have demonstrated how the non-invertible equivalence mappings presented in Section \ref{sec:noninvmap} are obtained through the symmetry-based method.

\section{Concluding remarks}
In this paper, by taking into account an application of the symmetry-based method \cite{blumanrafamsml2021a,blumanrafamsml2021b} to linear PDE systems introduced in \cite{blumanzuhal}, we have determined new non-invertible mappings for linear wave equations with variable wave speeds. As illustrative results: (1) We presented a new equivalence mapping of a linear wave equation with a variable wave speed nonlocally related to the linear wave equation (\ref{waveeqcx}) for any wave speed $c(x)$ through the symmetry-based method, extending results obtained in \cite{blumanrafamsml2021a,blumanrafamsml2021b}. (2) We determined non-invertible mappings of the linear wave equation (\ref{waveeq}) to the linear wave equation with a constant wave speed (\ref{ctecoeffpde}) for a new class of variable wave speeds $c(x,t)$. Hence, we obtained the general solutions of these linear wave equations with the new class of variable wave speeds $c(x,t)$. 

\section*{Acknowledgments}
We thank the Natural Sciences and Engineering Research Council of Canada for financial support.


\begin{thebibliography}{xx}

\bibitem{ancobluman1997}
S.C. Anco, G.W. Bluman, Nonlocal symmetries and nonlocal conservation laws of Maxwell's equations, J. Math. Phys. \textbf{38}, 3508--3532, 1997.

\bibitem{ancothe2005}
S.C. Anco, D. The, Symmetries, conservation laws, and cohomology
of Maxwell's equations using potentials, Acta Appl. Math. \textbf{89}, 1--52, 2005.

\bibitem{blok2021}
D.I. Blokhintsev, \textit{Acoustics of a Nonhomogeneous Moving Medium}, translated from the English by the NACA, TM 1399, 1956. Originally published as \textit{Akustika Neodnorodnoi Dvizhushcheisya Sredy, Ogiz, Gosudarstvennoe Izdatel'stvo, Tekhniko-Teoreticheskoi Literatury, Moskva, Leningrad, 1946}.

\bibitem{bluman1983}
G.W. Bluman, On mapping linear partial differential equations to constant coefficient equations, SIAM J. Appl. Math. \textbf{43(6)}, 1259--1273, 1983.

\bibitem{bluman2007}
G.W. Bluman, A.F. Cheviakov, Nonlocally related systems, linearization and nonlocal symmetries for the nonlinear wave equation, J. Math. Anal. Appl. \textbf{333}, 93--111, 
2007.

\bibitem{blumanchevianco}
G.W. Bluman, A.F. Cheviakov, S.C. Anco, \textit{Applications of Symmetry Methods to Partial Differential Equations}, Springer, New York, NY, 2010.

\bibitem{bluman2006}
G.W. Bluman, A.F. Cheviakov, N.M. Ivanova, Framework for nonlocally related partial differential equation systems and nonlocal symmetries: Extension, simplification, and examples, J. Math. Phys. \textbf{47}, 113505, 2006.

\bibitem{blumanrafa2025}
G.W. Bluman, R. de la Rosa, The natural extension to PDEs of Lie's reduction of order algorithm for ODEs, Commun. Nonlinear Sci. Numer. Simulat. \textbf{140}, 108438, 2025.

\bibitem{blumanrafamsml2021a}
G.W. Bluman, R. de la Rosa, M.S. Bruz\'on, M.L. Gandarias, Differential invariant method for seeking nonlocally related systems and nonlocal symmetries. I: General theory and examples, Proc. R. Soc. A \textbf{477}, 20200908, 2021.

\bibitem{blumanrafamsml2021b}
G.W. Bluman, R. de la Rosa, M.S. Bruz\'on, M.L. Gandarias, Differential invariant method for seeking nonlocally related systems and nonlocal symmetries. II: Connections with the conservation law method, Proc. R. Soc. A \textbf{477}, 20200909, 2021.

\bibitem{blumankumei1987}
G.W. Bluman, S. Kumei, On invariance properties of the wave equation, J. Math. Phys. \textbf{28}, 307--318, 1987.

\bibitem{blumankumei}
G.W. Bluman, S. Kumei, \textit{Symmetries and Differential Equations}, Springer-Verlag, New York, NY, 1989.

\bibitem{blumankumeireid1988}
G.W. Bluman, S. Kumei, G.J. Reid, New classes of symmetries for partial differential equations, J. Math. Phys. \textbf{29}, 806--811, 1988.

%
\bibitem{blumanyang2013}
G.W. Bluman, Z. Yang, A symmetry-based method for constructing nonlocally related partial differential equation systems, J. Math. Phys. \textbf{54}, 093504, 2013.


\bibitem{blumanzuhal}
G.W. Bluman, Z.K. Y\"{u}zba\c{s}i, How symmetries yield non-invertible mappings of linear partial differential equations, J. Math. Anal. Appl. \textbf{491(2)}, 124354, 2020.

\bibitem{brek1976}
L.M. Brekhovskikh, \textit{Waves in Layered Media}, Academic Press, New York, NY, 1980.

\bibitem{kaptsov2024}
O.V. Kaptsov, Representation of linear waves in inhomogeneous media, 2024. [https://doi.org/10.48550/arXiv.2410.13197]

\bibitem{mei2005}
C.C. Mei, M.A. Stiassnie, D.K.P. Yue, \textit{Theory and Applications of Ocean Surface Waves}, World Scientific, 2005.

\bibitem{melnikov2024}
I.E. Melnikov, E.N. Pelinovsky, Linear waves on shallow water slowing down near the shore over uneven bottom, Fluid Dyn. \textbf{59}, 260-269, 2024.

\bibitem{pelinovsky2022}
E.N. Pelinovsky, O.V. Kaptsov, Traveling waves in shallow seas of variable depths, Symmetry \textbf{14}, 1448, 2022.

\bibitem{seymour1987}
B. Seymour, E. Varley, Exact representations for acoustical waves when the sound speed varies in space and time, Stud. Appl. Math. \textbf{76}, 1-35, 1987.

\bibitem{yang2013} 
Z. Yang, Nonlocally Related Partial Differential Equation Systems, the Nonclassical Method and Applications,  Ph.D. Thesis, University of British Columbia, 2013.

%
%
%
%
%
%
%
%
%
%
%
%
%
%
%
\end{thebibliography}
\end{document}